\documentclass[11pt]{article}
\usepackage{amsfonts}
\usepackage{amssymb}
\usepackage{fullpage}
\usepackage{xy}
\xyoption{all}

\title{\bf{On the Ricci curvature of normal metrics on Biquotients}}
\author{Lorenz J. Schwachh\"ofer\footnote{Research supported by the 
Schwerpunktprogramm Differentialgeometrie of the Deutsche
Forschungsgesellschaft}}

\date{February 21, 2007}
\begin{document}
\maketitle


\newtheorem{thm}{Theorem}[section]
\newtheorem{lem}[thm]{Lemma}
\newtheorem{prop}[thm]{Proposition}
\newtheorem{df}[thm]{Definition}
\newtheorem{cor}[thm]{Corollary}
\newtheorem{rem}[thm]{Remark}
\newtheorem{ex}[thm]{Example}
\newenvironment{proof}{\medskip
\noindent {\bf Proof.}}{\hfill \rule{.5em}{1em}\mbox{}\bigskip}

\def\eps{\varepsilon}

\def\GlR#1{\mbox{\it Gl}(#1,\R)}
\def\glR#1{{\frak gl}(#1,\R)}
\def\GlC#1{\mbox{\it Gl}(#1,\C)}
\def\glC#1{{\frak gl}(#1,\C)}
\def\Gl#1{\mbox{\it Gl}(#1)}

\def\ov{\overline}
\def\ot{\otimes}
\def\und{\underline}
\def\w{\wedge}
\def\ra{\rightarrow}
\def\lra{\longrightarrow}
\def\less{\prec}

\newcommand{\al}{\alpha}
\newcommand{\be}{\beta}
\newcommand{\ga}{\gamma}
\newcommand{\la}{\lambda}
\newcommand{\om}{\omega}
\newcommand{\Om}{\Omega}
\renewcommand{\th}{\theta}
\newcommand{\Th}{\Theta}

\def\pair#1#2{Q\left(#1,#2\right)}
\def\big#1{\displaystyle{#1}}

\def\N{{\Bbb N}}
\def\Z{{\Bbb Z}}
\def\R{{\Bbb R}}
\def\C{{\Bbb C}}
\def\CP{{\Bbb C} {\Bbb P}}
\def\HP{{\Bbb H} {\Bbb P}}
\def\P{{\Bbb P}}
\def\Q{{\Bbb Q}}

\def\O{{\cal O}}
\def\H{{\cal H}}
\def\V{{\cal V}}
\def\F{{\cal F}}

\def\so{{\frak {so}}}
\def\co{{\frak {co}}}
\def\su{{\frak {su}}}
\def\uu{{\frak {u}}}
\def\sl{{\frak {sl}}}
\def\sp{{\frak {sp}}}
\def\csp{{\frak {csp}}}
\def\spin{{\frak {spin}}}
\def\g{{\frak g}}
\def\h{{\frak h}}
\def\k{{\frak k}}
\def\m{{\frak m}}
\def\n{{\frak n}}
\def\t{{\frak t}}
\def\s{{\frak s}}
\def\z{{\frak z}}
\def\p{{\frak p}}
\def\L{{\frak L}}
\renewcommand{\l}{{\frak l}}
\def\X{{\frak X}}
\def\gl{{\frak {gl}}}
\def\hol{{\frak {hol}}}
\renewcommand{\frak}{\mathfrak}
\renewcommand{\Bbb}{\mathbb}

\def\be{\begin{equation}}
\def\ee{\end{equation}}
\def\bi{\begin{enumerate}}
\def\ei{\end{enumerate}}
\def\ba{\begin{array}}
\def\ea{\end{array}}
\def\bea{\begin{eqnarray}}
\def\eea{\end{eqnarray}}
\def\ben{\begin{enumerate}}
\def\een{\end{enumerate}}

\def\dq{\slash \!\!\!\! \slash}
\def\codim{\mbox{\rm codim}}

\begin{abstract}
\noindent
We show that any normal metric on a closed biquotient with finite fundamental 
group has positive Ricci curvature.
\end{abstract}

\section{Introduction}

One of the classical problems in differential geometry is the investigation 
of closed manifolds which admit Riemannian metrics with given lower 
bounds for the sectional or the Ricci curvature. As a sample question, 
one may ask for manifolds with a metric of both nonnegative sectional and 
positive Ricci curvature. Note that such a manifold has necessarily a finite 
fundamental group by the Bonnet-Myers theorem.

By a recent result of B\"ohm and Wilking (\cite{BW}), a closed manifold 
with finite fundamental group which carries a metric of nonnegative 
sectional curvature also admits metrics of positive Ricci curvature. These 
metrics are obtained by a deformation of the original metric using the Ricci 
flow, so we may assume in addition that these metrics have {\em almost} 
nonnegative sectional curvature.

In general, however, manifolds which have nonnegative sectional and 
positive Ricci curvature simultaneously seem to be hard to come by. One 
class of such manifolds are compact homogeneous spaces $G/H$ with a 
{\em normal} metric, i.e., a metric for which the canonical submersion $G \ra 
G/H$ is Riemannian where $G$ is equipped with a biinvariant metric. Such 
a metric always has nonnegative sectional curvature by O'Neill's formula, and 
it is not hard to see that it has positive Ricci curvature as long as $G/H$ has 
finite fundamental group (\cite{Na}).

Another class of manifolds of nonnegative sectional curvature are closed 
biquotients $G \dq H$ for  $H \subset G \times G$ with a {\em normal} metric, 
by which we mean - in analogy to the homogeneous case - a metric for which 
the canonical submersion $G \ra G \dq H$ is Riemannian, $G$ being equipped 
with a biinvariant metric. Normal metrics are in some sense the most canonical 
metrics on closed biquotients, so the question under which circumstances such 
a metric has positive Ricci curvature is very natural.

Surprisingly, not much was known about this question up to now. In \cite{Stud}, 
this positivity was shown under very restrictive circumstances; in  \cite{ST}, we 
showed that a normal metric on a biquotient with finite fundamental group has 
positive Ricci curvature on a dense open subset. The aim of the present article 
is to resolve this question completely.

\

\noindent {\bf Main Theorem} {\em Let $M := G \dq H$ be a biquotient
of a compact Lie group $G$, equipped with a normal metric. 
If $M$ has finite fundamental group then $M$ has positive Ricci curvature.}

\

Thus, any closed biquotient with finite fundamental group carries a metric of 
both nonnegative sectional and positive Ricci curvature. In general, we show 
that a closed biquotient with a normal metric is finitely isometrically covered by 
$M' \times T^k$, where $M'$ is a biquotient with a normal metric of positive Ricci 
curvature and $T^k$ is a flat torus. 

It is my pleasure to 
thank the Max-Planck-Institute 
for Mathematics in the Sciences in Leipzig and 
the Unversit\'e Libre de Bruxelles for 
their hospitality during the preparation
of parts of this work.

\section{Normal metrics on Biquotients} \label{sec:normal-biquotients}

Let $G$ be a compact connected Lie group with Lie algebra $\g$, and let 
$H \subset G \times G$ be a closed subgroup for which the action of $H$ on 
$G$ given by
\[
(h_1, h_2) \cdot g := h_1 g h_2^{-1}
\]
is free and hence the quotient $G \dq H$ is a manifold, called the {\em 
biquotient of $G$ by $H$}.

Let $Q$ be an $Ad_G$-invariant inner product on $\g$ which induces a 
biinvariant Riemannian metric of nonnegative sectional curvature on $G$.
Since the action of $H$ on $G$ preserves $Q$, there is a unique Riemannian 
metric $g_Q$ on $G \dq H$ for which the projection
\[
\pi: G \longrightarrow G \dq H
\] 
becomes a Riemannian submersion. Evidently, $g_Q$ has nonnegative 
sectional curvature by O'Neill's formula. 

\begin{df}
A metric $g_Q$ on $G \dq H$ induced by the $Ad_G$-invariant inner product $Q$ 
on $\g$ as above is called normal.
\end{df}

If $H = L \times K \subset G \times G$ then $H$ acts freely on $G$ iff $Ad_g 
L \cap K = \{e\}$ for all $g \in G$. In this case, we denote the biquotient $G 
\dq H$ also by $L \backslash G/K$. In fact, {\em any} biquotient can be 
considered to be of this form with $G/K$ being a {\em symmetric space}. 
Namely, there is a canonical diffeomorphism
\be \label{eq:general-biquotient}
G \dq H \longleftrightarrow H \backslash (G \times G)/\Delta G.
\ee
If we consider the normal metric $g_{Q \oplus Q}$ on $H \backslash (G \times 
G)/\Delta G$, then one sees easily that the diffeomorphism in 
(\ref{eq:general-biquotient}) becomes an isometry. Therefore, since 
$(G \times G)/\Delta G$ is a symmetric space, we shall from now on always 
consider biquotients of the form
\be \label{eq:general-restriction}
M = H \backslash G/K,\ \ \ \ \mbox{where $G/K$ is a symmetric space},
\ee
and we have the corresponding involution $\tau$, given by 
\be \label{eq:def-tau}
\tau: \g \ra \g,\ \ \ \tau|_\k = Id_\k,\ \ \ \ \tau_{\k^\perp} = -Id_{\k^\perp}.
\ee 
For $g \in G$, we define the $Q$-orthogonal decomposition
\be \label{eq:decomposition}
\ba{lll}
\g = \V_g \oplus \H_g, & \mbox{where} & \V_g := Ad_{g^{-1}}(\h) \oplus \k,
\ea
\ee
with the Lie algebras $\h$ and $\k$ of $H \subset G$ and $K \subset G$,
respectively. Note that $dL_g(\V_g) = \ker(d\pi_g)$, so that $T_gG = dL_g(\V_g) 
\oplus dL_g(\H_g)$ is the canonical decomposition of $TG$ into the vertical 
and horizontal subspace w.r.t. the Riemannian submersion $\pi: G \ra G\dq H$. 
In particular, the dimensions of $\V_g$ and $\H_g$ are independent of $g$.

\section{Biquotients with a flat point}  \label{sec:flat-points}

As a first step towards our main result, we investigate biquotients with flat
points, i.e., points where the curvature vanishes entirely. Since $\pi: G \ra 
G \dq H$ is a Riemannian submersion, hence $Sec^{G \dq H}_{\pi(g)}(\und v \w
\und w) \geq Sec^G_g(v \w w) = \frac14 ||[v, w]||^2$, where $v, w \in 
\H_g$ are horizontal lifts of $\und v$ and $\und w$, respectively, it follws 
that $\pi(g) \in G \dq H$ is a flat point only if $[\H_g, \H_g] = 0$. In this 
case, we can prove the following

\begin{thm} \label{thm:flat-points}
Let $G \dq H$ be a biquotient with a normal metric. If $G \dq H$ has a flat 
point, i.e., if there is a $g \in G$ with $[\H_g, \H_g] = 0$, then $\H_g 
\subset \z(\g)$. In this case, $G \dq H$ is finitely covered by a flat torus.
\end{thm}

\begin{proof} We assume w.l.o.g. that $g = e$. After passing to a 
finite cover, we further assume that $G = G_s \times Z$, where $G_s$ is 
semi-simple and simply connected and $Z = Z(G)$ is the identity component of the 
center of $G$. Moreover, we may assume that $H \subset G \times G$ is connected. 
We have an orthogonal splitting of the Lie algebra $\g = \g_s \oplus \z(\g)$, and we 
let $T \subset \ov{\exp \H_e} \subset G$ be a torus with Lie algebra $\t$ such that 
$\g = \V_e \oplus \t$ as a vector space. In particular, $\dim T = \dim \H_e$. (We use 
$T$ since the abelian subgroup $\exp \H_e \subset G$ may be not closed). Now 
consider the map $p: T \ra G \dq H$ which makes the following diagram commute.
\be \label{diagram2}
\xymatrix{
& G \ar[dd]^{\mbox{\footnotesize $\pi$}}\\
T \ar@{^{(}->}[ur] \ar[dr]^{p} & \\ & G \dq H
}
\ee

\noindent Let $g \in T$. Since $T \subset \ov{\exp(\H_e)}$ and the latter is abelian, we
have
\[
Q(Ad_{g^{-1}}(\h), \H_e) = Q(\h, Ad_g(\H_e)) = Q(\h, \H_e) = 0,\ \ \ \
\mbox{and} \ \ \ \ Q(\k, \H_e) = 0.
\]
Thus, $Q(\V_g, \H_e) = 0$, which implies that $\V_g = \V_e$ for all $g \in
T$. In particular, $\t$ is transversal to $\V_g$ for all $g \in T$. As $\V_g$ 
is by definition the kernel of the differential $d\pi_g: T_gG \ra
T_{p(g)}G \dq H$ for $g \in T$, it follows that $dp_g: dL_g(\t) \hookrightarrow 
T_gG \ra d_{p(g)}G \dq H$ is an isomorphism, hence $p: T \ra G \dq H$ is a 
finite covering map, hence the induced map $p_*: \pi_1(T) \ra \pi_1(G \dq H)$ 
is injective with index $[\pi_1(G \dq H) : p_*(\pi_1(T))] < \infty$.

Since $G = G_s \times Z$ with $G_s$ simply connected, it follows that the
projections $pr_Z: G \ra Z$ induces an isomorphism of fundamental groups, 
hence from the homotopy exact sequence
\[
\xymatrix{
& \pi_1(T) \ar[d]_{\mbox{\footnotesize $(pr_Z)_*$}}
\ar[dr]^{\mbox{\footnotesize $p_*$}} & & \\
\pi_1(H) \ar[r]_{\mbox{\footnotesize $(pr_Z)_*$}} & \pi_1(Z)
\ar[r]^{\mbox{\footnotesize $\pi_*$}} & \pi_1(G \dq H) \ar[r] & 0
}
\]
and from $[\pi_1(G \dq H) : p_*(\pi_1(T))] < \infty$ we conclude that $\dim pr_{\z(\g)}(\t) 
= \dim \t$ and that $\z(\g) = pr_{\z(\g)}(\V_e) \oplus pr_{\z(\g)}(\t)$. Then, for 
dimensional reasons, it follows that $\V_e = pr_{\z(\g)}(\V_e) \oplus \g_s$ and therefore, 
$\H_e \subset \g_s^\perp = \z(\g)$.
\end{proof}

\section{Biquotients with Ricci-flat directions}  \label{sec:Ricci-flat}

Let us suppose that at a point $p \in G \dq H$ the Ricci curvature is not positive, 
i.e., there are tangent vectors $0 \neq \und v \in T_pG \dq H$ with $Sec^{G \dq 
H}_{\pi(g)}(\und v \w \und w) = 0$ for all $\und w \in T_pG\dq H$. W.l.o.g. we assume 
that $p = \pi(e)$. As in the preceding section, O'Neill's formula then implies that $[v, 
\H_e] = 0$ where $v \in \H_e$ is the horizontal lift of $\und v$. To simplify our notation 
we shall from now on omit the subscript and thus let $\V := \V_e = \h \oplus \k$ and 
$\H := \H_e = \V^\perp$. Thus, we define the space of Ricci flat directions at $\pi(e)$
\be \label{eq:define-F}
\F := \{ v \in \H \mid [v, \H] = 0 \} \subset \H
\ee
and assume that $\F \neq 0$.

\begin{prop} \label{prop:MainThm}
Let $M := H \backslash G/K$ be a normal biquotient as in (\ref{eq:general-restriction}), 
and suppose that $\F \subset \z(\g)$. Then $M$ is finitely isometrically covered by 
$M' \times T^k$, where $M'$ is a biquotient with finite fundamental group and a normal 
metric, and $T^k$ is a flat torus of dimension $k = \dim \F$.
\end{prop}

\begin{proof} 
After replacing $G$ by a finite cover if necessary, we may assume that $G = G_s \times 
Z$ where $G_s$ is semisimple and $Z$ is the identity component of the center with Lie 
algebra $\z(\g)$. Moreover, after replacing $M$ by a finite isometric cover, we may 
assume that $H, K \subset G$ are connected. Let $pr_Z : G = G_s \times Z \ra Z$ be the 
canonical projection.

Note that the map $\jmath: H \times K \ra Z$, $(h,k) \mapsto pr_Z(h) \cdot pr_Z(k)$ is a 
homomorphism, hence its image $\jmath(H \times K) \subset Z$ is a torus. Also, by 
hypothesis, it follows that $\F = \z(\g) \cap \H$, so that $\z(\g) = pr_{\z(\g)}(\h \oplus \k) 
\oplus \F$.

Let $G' := G_s \times \jmath(H \times K)$. Then evidently, $G' \lhd G$, and $G \cong 
G' \times T^k$, where $T^k \subset Z$ is a torus with Lie algebra $\t \subset \z(\g)$ of 
dimension $k = \dim(\F)$.

Choose a linear map $p: \g \ra \g$ such that $p|_{\g'} = Id_{\g'}$ and $p(\t) = \F$, and 
define the inner product $\tilde Q := p^*(Q)$. Since $\tilde Q|_{\g'} = Q|_{\g'}$ and $\tilde 
Q(\g_s, \z(\g)) = 0$, it follows that $\tilde Q$ is biinvariant as well. Note that $H \times K 
\subset G'$, hence we have a canonical isometry
\[
(H \backslash G/K, g_Q) \longleftrightarrow ((H \backslash G'/K) \times T^k, g_{\tilde Q}),
\]
where the latter is a product metric, and the metric on the torus factor is flat. Also, by 
construction, $M' := ((H \backslash G'/K), g_{\tilde Q}) = ((H \backslash G'/K), g_Q)$ is a 
biquotient with a normal metric which has positive Ricci curvature at $\pi(e)$. Therefore, by 
the deformation result of Ehrlich \cite{E}, $M'$ admits a metric of positive Ricci curvature 
and hence has finite fundamental group.
\end{proof}

\begin{df} \label{def:aberrant}
Let $M := H \backslash G/K$ be a normal biquotient as in (\ref {eq:general-restriction}).
We say that $M$ has {\em non-central flats} if  $\F \not\subset \z(\g)$. If $M$ has 
non-central flats, then we call $M$ {\em minimal} if for any normal biquotient $M' = 
H' \backslash G'/K'$ with non-central flats with $\dim M' \leq \dim M$ and $\dim G' \leq 
\dim G$ we have $\dim M' = \dim M$ and $\dim G' = \dim G$.
\end{df}

Thus, our main theorem will follow if we can show that there do not exist normal biquotients 
with non-central flats which shall be our goal for the remainder of this paper and will be 
achieved in Theorem~\ref{thm:MainTheorem}. Indeed, once this is shown it will follow from 
Proposition~\ref{prop:MainThm} that a biquotient with a normal metric with Ricci flat directions 
has infinite fundamental group.

\begin{prop} \label{prop:irreducible}
Suppose that $M = H \backslash G/K$ is a minimal biquotient with non-central flats. Then $G/K$ is a 
compact irreducible symmetric space.
\end{prop}

\begin{proof} Suppose that $\g = \g_1 \oplus \g_2$ is a $\tau$-invariant decomposition 
where $\tau: \g \ra \g$ denotes the symmetric involution from (\ref{eq:def-tau}). After 
passing to a covering if necessary, we may assume w.l.o.g. that $G/K = G_1/K_1 \times 
G_2/K_2$ where $G = G_1 \times G_2$ and $K_i := K \cap G_i$. Then $\k = (\k \cap \g_1) 
\oplus (\k \cap \g_2) =: \k_1 \oplus \k_2$. Let $pr_i: \g \ra \g_i$ be the canonical projection. 
Let $\H_i := pr_i(\H)$ and $\F_i := pr_i(\F)\subset \H_i \subset \k_i^\perp$. If $\F_i \subset 
\z(\g_i)$ for $i = 1,2$ then we would have $\F \subset \z(\g)$ which we assumed {\em not} to be the case, 
thus we may assume w.l.o.g. that $\F_1 \not\subset \z(\g_1)$.

Now consider the biquotient $M' := (H \cap G_1) \backslash G_1/K_1$. Its horizontal space 
is $\k_1^\perp \cap (\h \cap \g_1)^\perp = \H_1$, hence $\F_1 \subset \H_1$ is the space 
of flat directions of $M'$. Since $\F_1 \not\subset \z(\g_1)$, it follows that $M'$ has non-central flats, 
and $\dim M' = \dim \H_1 \leq \dim \H = \dim M$. But $\dim G_1 \leq \dim G$, hence the 
minimality of $M$ implies that $\dim G_1 = \dim G$, i.e., $G_2 = 1$.
\end{proof}

The following simple facts will be useful later on.

\begin{lem}\label{lem:generatedideal}
Let $\g$ be a Lie algebra with a biinvariant inner product $Q$, and let $L \subset \g$ be a 
linear subspace. Let $L_1 := [L, \g]$. Then the ideal generated by $L$ is the linear span of 
$L$, $L_1$ and $[L_1, L_1]$.
\end{lem}

\begin{proof}
Since $ad_x: \g \ra \g$ is skew symmetric for all $x \in \g$, it follows that the images of 
$ad_x$ and $(ad_x)^2$ are equal, hence $[L, L_1] = L_1$.

Let $I := span(L, L_1, [L_1,L_1]) \subset \g$. Then $Q([I^\perp,L], \g) = Q(I^\perp, [L,\g]) = 
Q(I^\perp,L_1) = 0$, so that $[I^\perp,L] = 0$. Next, $[I^\perp, L_1] = [I^\perp, [L, L_1]] = 
[L, [I^\perp, L_1]] \subset [L, \g] = L_1$. On the other hand, $Q([I^\perp, L_1], L_1) = 
Q(I^\perp, [L_1,L_1]) = 0$, hence $[I^\perp, L_1] = 0$. Then the Jacobi identity implies 
that $[I^\perp, [L_1,L_1]] = 0$ as well.

Thus, $[I, I^\perp] = 0$, hence $Q([I, \g], I^\perp) = Q(\g, [I, I^\perp]) = 0$, which shows 
that $I \lhd \g$ is an ideal. On the other hand, any ideal containing $L$ must also contain 
$I$ which shows the claim.
\end{proof}

\begin{lem} \label{lem:composition=0}
Let $G/K$ be a compact irreducible symmetric space with Lie algebra $\g$. Let $\s 
\subsetneq \g$ be a proper subalgebra, and let $\l_1, \l_2 \subset \s$ be subalgebras 
such that $(ad_{\l_2} \circ ad_{\l_1})|_{\s^\perp}$ vanishes. If there is an $ad_\s$-invariant 
and $\tau$-invariant subspace $0 \neq W \subset \s^\perp$ such that $[\l_1, W] = W$, 
then $\l_2 = 0$.
\end{lem}

\begin{proof}
First observe that $[\l_2, W] = [\l_2, [\l_1, W]] \subset (ad_{\l_2} \circ ad_{\l_1})(\s^\perp) = 
0$, i.e., $[\l_2, W] = 0$.

Next, we wish to show that $[\l_2, [W, \g]] = 0$. For this, we observe that $[\l_2, [W, \s]] = 
[\l_2, W] = 0$, and $[\l_2, [W, W]] = 0$ by the Jacobi identity. Thus, it remains to show that 
$[\l_2, [W, W']] = 0$ where $W' := W^\perp \cap \s^\perp$. Since $Q([W, W'], \s) = Q(W', 
[\s, W]) \subset Q(W', W) = 0$, it follows that $[W, W'] \subset \s^\perp$. Thus,
\[
\ba{lll}
[\l_2, [W, W']] & = & [\l_2, [[\l_1, W], W']]\\
& \subset & \underbrace{[\l_2, [\l_1, [W, W']]]}_{\subset (ad_{\l_2} \circ ad_{\l_1})(\s^\perp)  
= 0} + [\l_2, [W, [\l_1,W']]\\
& \subset & [\underbrace{[\l_2, W]}_{=0}, [\l_1,W']] + [W, \underbrace{[\l_2, [\l_1,W']]}_{\subset 
(ad_{\l_2} \circ ad_{\l_1})(\s^\perp)  = 0}] = 0.
\ea
\]

Therefore, $[\l_2, W] = [\l_2, [W, \g]] = 0$, and hence, according to 
Lemma~\ref{lem:generatedideal}, $[\l_2, I] = 0$ where $I \lhd \g$ is the ideal generated by 
$W$. Since $W \neq 0$ is $\tau$-invariant, so is $I \neq 0$, hence the irreducibility of $G/K$ 
implies that $I = \g$, i.e., $\l_2 \subset \z(\g) = 0$.
\end{proof}

Suppose that $M = H \backslash G/K$ is a minimal biquotient with non-central flats, so that $G/K$ is an 
irreducible symmetric space by Proposition~\ref{prop:irreducible}, and $\F \neq 0$. Let
\[
S := Stab(\F) = \{ g \in G \mid Ad_g|_\F = Id_\F \} \subsetneq G.
\]
Evidently, $S$ is compact and $\tau$-invariant, and the Lie algebra $\s$ of $S$ is the centralizer 
of $\F$, i.e.,
\[
\s = \{ v \in \g \mid [v, \F] = 0 \} \subsetneq \g.
\]

Then $\H \subset \s$ by (\ref{eq:define-F}), and we let $\l \subset \s$ be the Lie subalgebra 
generated by $\H$. Since $\H$ is $\tau$-invariant w.r.t. the involution $\tau: \g \ra \g$ from 
(\ref{eq:def-tau}), so is $\l$. Evidently, $\F = \z(\l)$, hence $\l = \F \oplus \l_s$ where $\l_s$ is 
semi-simple, and $(\l_s, \tau)$ is a symmetric pair without Euclidean factor.

If $\l_s = 0$, i.e., $\H = \F$ is abelian, then $\F \subset \z(\g)$ by Theorem~\ref{thm:flat-points}. 
Thus, for biquotients with non-central flats, we may assume that $\l_s \neq 0$.

\begin{prop} \label{prop:l_s preserves s^perp}
Let $H \backslash G/K$ be a minimal biquotient with non-central flats. Let $\l_s = \l_1 \oplus \l_2$ be a decomposition into $\tau$-invariant ideals such that
\be \label{eq:decompose-H}
\H = \F \oplus (\H \cap \l_1) \oplus (\H \cap \l_2).
\ee
Then either $\l_k = 0$ or $[\l_k, \s^\perp] = \s^\perp$ for $k = 1,2$. In particular, $[\l_s, \s^\perp] 
= \s^\perp$.
\end{prop}

\begin{proof}
Let us assume that, say, $\l_1 \neq 0$. Since $\l_1 \subset \l \subset \s$, it is evident that 
$[\l_1, \s^\perp] \subset \s^\perp$. 

Let $G' := \{ g \in G \mid Ad_g(\l_1 \cap \k^\perp) = \l_1 \cap \k^\perp\}$. 
Evidently, $G' \subset G$ is a closed subgroup, hence also compact. Since 
$G'$ and hence its Lie algebra $\g'$ are $\tau$-invariant, we have the 
decomposition $\g' = (\g' \cap \k) \oplus (\g' \cap \k^\perp) =: \g_0'
\oplus \g_1'$, where
\[
\g_0' = \{ k \in \k \mid [k, \l_1] \subset \l_1\},\ 
\ \ \ \ \mbox{and}\ \ \ \ \
\g_1' = \{ x \in \k^\perp \mid [x, \l_1] = 0\},
\]
as $\l_1$ is generated by $\l_1 \cap \k^\perp$. Since $\l_1$ is semi-simple, 
we have
\be \label{eq:g1'-perp}
Q(\g_1', \l_1) = Q(\g_1', [\l_1, \l_1]) = Q(\underbrace{[\g_1', \l_1]}_{=0}, \l_1) = 0.
\ee
Since $\F \oplus \l_2 \subset \g'$, we conclude by (\ref{eq:decompose-H}) and 
(\ref{eq:g1'-perp}) that
\be \label{eq:decompose-g'1}
\g' = \F \oplus (\H \cap \l_2) \oplus (\g' \cap \H^\perp) = \F  \oplus (\H \cap \l_2) \oplus 
(\g' \cap (\h \oplus \k)).
\ee

Let $x \in pr_{\k^\perp} (\g' \cap (\h \oplus \k))$. Then $x \in \g_1'$ and 
there is a $k \in \k$ such that $x + k \in \h$. We have on the one hand 
$[x + k, \H \cap \l_1] \subset [\h, \H] \subset [\h, \h^\perp] \subset \h^\perp$. On the other 
hand, $[x, \H \cap \l_1] \subset [\g_1', \l_1] = 0$, hence $[x + k, \H \cap \l_1] = 
[k, \H \cap \l_1] \subset [\k, \k^\perp] \subset \k^\perp$.

Thus, $[x + k, \H \cap \l_1] = [k, \H \cap \l_1] \subset \h^\perp \cap \k^\perp = \H$, and 
$Q([k, \H \cap \l_1], \F \oplus \l_2) = Q(k, [\H \cap \l_1, \F \oplus \l_2]) = 0$, i.e., $[k, \H \cap \l_1] 
\subset \H \cap \l_1$. Since $\H \cap \l_1$ generates $\l_1$, it follows that $[k, \l_1] \subset 
\l_1$ and hence, $k \in \g_0'$, i.e., $x + k \in \h \cap \g'$.

Therefore, $pr_{\k^\perp}(\g' \cap (\h \oplus \k)) = pr_{\k^\perp}(\g' \cap \h)$,
so that we can refine (\ref{eq:decompose-g'1}) to the (not orthogonal) direct sum 
decomposition
\[
\g' = \F \oplus (\H \cap \l_2) \oplus (\g' \cap \h) \oplus (\g' \cap \k).
\]
Thus, for the biquotient $(G' \cap H) \backslash G' / (G' \cap K)$, the
horizontal space at $e \in G'$ is $\F \oplus (\H \cap \l_2) \subsetneq \H$ as $\l_1 \neq 0$, 
i.e., $\dim (G' \cap H) \backslash G' / (G' \cap K) < \dim H \backslash G/K$, and $\dim 
G' \leq \dim G$. Therefore, the minimality of $H \backslash G/K$ implies that $(G' \cap H) 
\backslash G' / (G' \cap K)$ cannot have non-central flats. Since the centralizer of $\F \oplus (\H \cap 
\l_2)$ equals $\F$, this implies that $\F \subset \z(\g')$ or, equivalently, $\g' \subset \s$.

Let $v \in \s^\perp \cap [\l_1, \s^\perp]^\perp$. Then
\[
Q(\underbrace{[v, \l_1]}_{\subset \s^\perp}, \s^\perp) = Q(v, [\l_1, \s^\perp]) = 0,
\]
hence $[v, \l_1] = 0$, i.e., $v \in \g' \subset \s$. This implies that $v = 0$, i.e. 
$[\l_1, \s^\perp] = \s^\perp$, which shows our assertion.
\end{proof}

Let us now investigate the structure of the Lie algebra $\g = \s \oplus
\s^\perp$ in some more detail. We fix once and for all a {\em generic} 
element $v_0 \in \F$, meaning that $[x, \F] = 0$ iff $[x, v_0] = 0$. 
Thus, the skew-symmetric map $ad_{v_0}|_{\s^\perp}: \s^\perp \ra \s^\perp$ 
is invertible, hence there is a unique orthogonal complex structure $J: 
\s^\perp \ra \s^\perp$ and a positive definite symmetric map $P_0: \s^\perp 
\ra \s^\perp$ such that 
\be \label{eq:define P_0}
ad_{v_0}|_{\s^\perp} = J P_0: \s^\perp \longrightarrow \s^\perp.
\ee

\begin{lem} \label{lem:Lagrange}
Let $V_1 := pr_{\s^\perp} (\h)$ and $V_2 := pr_{\s^\perp}(\k)$. Then the 
following hold:
\bi
\item
$Q(\H, [V_i, V_i]) = 0$ for $ i = 1,2$, and hence $[\H, V_i] \subset V_i^\perp \cap \s^\perp$.
\item
$\s^\perp = V_1 \oplus V_2$ and $\dim V_1 = \dim V_2 = \frac12 \dim \s^\perp$.
\item
As a vector space, $\k = (\k \cap \s) \oplus V_2$.
\item
$P_0 V_2 = V_2$ and $V_2 \subset \s^\perp$ is {\em totally real w.r.t. $J$}, 
i.e. $J V_2 = V_2^\perp \cap \s^\perp$.
\item
$[\h \cap \s, V_1] \subset V_1$.
\ei
\end{lem}

\begin{proof}
Let $x_1,x_2 \in V_1$. Then there exist $k_1, k_2 \in \s \cap \k$ such that 
$x_i + k_i \in \h$. Thus, for $v \in \H$ it follows
\[
\ba{lll}
0 & = & Q(v, \underbrace{[x_1 + k_1, x_2 + k_2]}_{\in \h}) = Q(v, [x_1, x_2]) + 
Q(v, \underbrace{[x_1, k_2]}_{\in \s^\perp} + \underbrace{[k_1, x_2]}_{\in
\s^\perp} + \underbrace{[k_1, k_2]}_{\in \s \cap \k})\\ 
\\ & = & Q(v, [x_1, x_2]),
\ea
\]
which shows the first claim for $V_1$, and an analogous proof works for $V_2$. 

If $v_0 \in \F$ is generic, then the form $\om \in \Lambda^2(\s^\perp)^*$ 
defined by $\om(x,y) := Q(v_0, [x,y]) = Q([v_0,x], y)$ is non-degenerate; 
indeed, if $0 = \om(x,\s^\perp) = Q([v_0,x], \s^\perp)$ then $[v_0, x] = 0$ 
which implies that $x = 0$ as $v_0$ is generic.

Thus, by the first step, $V_i \subset \s^\perp$ is isotropic w.r.t. $\om$, hence
$\dim V_i \leq \frac12 \dim \s^\perp$. Moreover, if $x \in \s^\perp \cap 
V_1^\perp \cap V_2^\perp$, then $x \in \h^\perp \cap \k^\perp = \H \subset
\s$, so that $x = 0$. This shows that $V_1 + V_2 = \s^\perp$ and hence the 
second claim follows.

For the third assertion, suppose there is an $x \in V_2$ and $s_1, s_2 \in \s$ 
such that $x + s_1 \in \k$, $x + s_2 \in \k^\perp$. Then, since $v_0 \in \F 
\subset \k^\perp$ and $G/K$ is symmetric, it follows that $[v_0, x] = [v_0, 
x + s_i] \in \k^\perp \cap \k = 0$, hence $x = 0$. From this, it follows that
$\k = (\k \cap \s) \oplus V_2$.

Therefore, $[v_0,[v_0, V_2]] \subset [\k^\perp, [\k^\perp, \k]] \cap \s^\perp
\subset \k \cap \s^\perp = V_2$, and by (\ref{eq:define P_0}) this implies 
that $P_0^2 V_2 \subset V_2$. But $P_0$ is symmetric and invertible, hence 
$P_0 V_2 = V_2$. Thus, $J V_2 = J P_0 V_2 = [v_0, V_2] \subset [\k^\perp, \k]
\cap \s^\perp \subset \k^\perp \cap \s^\perp = V_2^\perp \cap \s^\perp$ which
shows the fourth part.

For the last assertion, let $h \in \h \cap \s$ and $x \in V_1$. Then there is 
an $s \in \s$ such that $x + s \in \h$ and hence $\h \ni [h, x + s] = [h,x] +
[h,s]$. But $[h,s] \subset \s$, whereas $[h,x] \in \s^\perp$, i.e., $[h,x] \in
pr_{\s^\perp}(\h) = V_1$. 
\end{proof}

By virtue of this lemma, we may regard $\s^\perp \cong V_2 \ot \C$. We 
denote the eigenspace decomposition of $\s^\perp$ w.r.t. $P_0$ by
\be \label{eq:decompose-V2}
\s^\perp = W_1 \oplus \ldots \oplus W_n,
\ee
with eigenvalues $0 < \la_1 < \ldots < \la_n$. Since $v_0 \in \z(\s) \cap \k^\perp$, 
it follows that $W_k$ is a complex  subspace which is invariant both under $\tau$ 
and $ad_\s$. Moreover, since $\s^\perp = V_1 \oplus V_2$, we have
\be \label{eq:define phi}
V_1 = \{ Jx + P_0^{-1} \phi x \mid x \in V_2\}
\ee
for some linear map $\phi: V_2 \ra V_2$. By abuse of notation we denote the complex 
linear extension  of $\phi$ also by $\phi: \s^\perp \ra \s^\perp$. Evidently, the latter map is 
$\tau$-invariant. 

\begin{lem} \label{properties phi}
The map $\phi: V_2 \ra V_2$ from (\ref{eq:define phi}) is symmetric, and $\phi$ and 
$P_0$ have no common eigenvectors.
\end{lem}

\begin{proof}
By Lemma~\ref{lem:Lagrange}.1, $Q([v_0, V_1], V_1) = 0$ and thus by 
(\ref{eq:define phi}), we have for $x,y, \in V_2$
\[
\ba{lll}
0 & = & Q(ad_{v_0}(Jx + P_0^{-1} \phi x), Jy + P_0^{-1} \phi y)\\
& = & Q(J P_0 (Jx + P_0^{-1} \phi x), Jy + P_0^{-1} \phi y)\\
& = & Q(- P_0 x + J \phi x, Jy + P_0^{-1} \phi y)\\
& = & -Q(P_0 x, P_0^{-1} \phi y) + Q(\phi x, y)\\
& = & -Q(x, \phi y) + Q(\phi x, y).
\ea
\]

For the second part, suppose that $0 \neq x \in V_2$ is such that $P_0 x = \la x$ and 
$\phi x = \mu x$. Then $V_1 \ni Jx + P_0^{-1} \phi x = (i + \mu/\la) x$, using the 
identification $\s^\perp \cong V_2 \ot \C$, and there is a $k \in \k \cap \s$ such that 
$(i + \mu/\la) x + k \in \h$.

Therefore, $Ad_{\exp t v_0} \h \ni e^{i t \la}(i + \mu/\la) x + k$. But for a suitable $t_0 \in 
\R$, we can achieve that $e^{i t_0 \la}(i + \mu/\la) \in \R$, i.e., 
$e^{i t_0 \la}(i + \mu/\la) x \in V_2$ and hence $0 \neq e^{i t_0 \la}(i + \mu/\la) x + k \in 
Ad_{\exp t_0 v_0} \h \cap \k$ which contradicts the biquotient property.
\end{proof}

Let $v \in \k^\perp \cap \s$. Then $[v, V_2] \subset [\s, \s^\perp] \cap 
[\k^\perp, \k] \subset \s^\perp \cap \k^\perp = JV_2$, where the last equation 
follows from Lemma~\ref{lem:Lagrange}.4. Therefore, there is a symmetric map 
$A_v: V_2 \ra V_2$ such that $[A_v, P_0] = 0$ and $ad_v|_{\s^\perp} = JA_v$.

\

\begin{lem} \label{lem:eigenvalues} Let $v \in \k^\perp \cap \s$ and let $A_v: V_2 \ra 
V_2$ be the symmetric linear map such that $ad_v|_{\s^\perp} = JA_v$.
\bi
\item If $v \in \H$ then $[P_0^{-1}A_v, \phi] = 0$.
\item If $v \in \V$ and $x \in \s^\perp$ is a common eigenvector of $P_0^{-1} A_v$ and
$\phi$, then $A_vx = 0$.
\ei
\end{lem}

\begin{proof} Let $v \in \H$ and $x, y \in V_2$, so that $J x + P_0^{-1} \phi 
x, J y + P_0^{-1} \phi y \in V_1$ by (\ref{eq:define phi}). By Lemma~\ref{lem:Lagrange}.1,
\[
\ba{llll}
0 & = & Q([v, J x + P_0^{-1} \phi x], J y + P_0^{-1} \phi y)\\
& = & Q(J A_v (J x + P_0^{-1} \phi x), J y + P_0^{-1} \phi y)\\
& = & Q(- A_v x + J A_v P_0^{-1} \phi x, J y + P_0^{-1} \phi y)\\
& = & - Q(A_v x, P_0^{-1} \phi y) + Q(A_v P_0^{-1} \phi x, y)\\
& = & -Q(\phi P_0^{-1} A_v x, y) + Q(P_0^{-1} A_v \phi x, y) & \mbox{since $\phi$
and $P_0$ are symmetric and $[A_v, P_0] = 0$}\\
& = & Q([P_0^{-1} A_v, \phi] x, y),
\ea
\]
which shows the first part. For the second, let $v \in \V \cap \k^\perp \cap \s$ and $0 \neq x
\in V_2$ such that $P_0^{-1} A_v x = \mu x$ and $\phi x = \nu x$. Let $k \in \k \cap \s$ be such 
that $v + k \in \h$. Thus, $J x + \nu P_0^{-1} x \in V_1$ by (\ref{eq:define phi}), and by
Lemma~\ref{lem:Lagrange}.5, $(ad_v + ad_k)(V_1) \subset V_1$, so that
\[
\ba{lll}
V_1 & \ni & (ad_v + ad_k)(Jx + \nu P_0^{-1}x)\\
& = & -A_v x + \nu J P_0^{-1} A_v x + J ad_k x + \nu P_0^{-1} ad_k x\\
& = & -\mu P_0 x + \nu \mu J x + J ad_k x + \nu P_0^{-1} ad_k x\\
& = & J(ad_k x + \nu \mu x) + P_0^{-1}(-\mu P_0^2 x + \nu\ ad_k x).
\ea
\]
Therefore, by (\ref{eq:define phi}), $-\mu P_0^2 x + \nu\ ad_k x = \phi(ad_k 
x + \nu \mu x)$ or 
\be \label{eq:lem-eigenvalues}
\nu\ ad_k x - \phi\ ad_k x = \mu (P_0^2 x + \nu \phi x) = \mu (P_0^2 + \nu^2) x.
\ee
Since $\phi$ is symmetric and $ad_k$ is skew-symmetric, we have
\[
Q(\nu\ ad_k x - \phi ad_k x, x) = \nu \underbrace{Q(ad_k x, x)}_{=0} - 
\underbrace{Q(ad_k x, \overbrace{\phi x}^{\nu x})}_{=0} = 0,
\]
and because of (\ref{eq:lem-eigenvalues}),
\[
\mu\ Q((P_0^2 + \nu^2) x, x) = 0.
\]
Since $P_0$ and hence $P_0^2 + \nu^2$ is positive definite, this implies that $\mu = 0$, 
i.e., $A_v x = 0$. 
\end{proof}

We call a subspace $U \subset \s^\perp$ {\em $(\tau, \phi, \l)$-invariant} if it is invariant 
under $\tau$, $\phi$ and $ad_\l$. In particular, $U$ must be $ad_{v_0}$-invariant, so that 
each non-zero $(\tau, \phi, \l)$-invariant subspace is complex and has the form
\[
U = U_1 \oplus \ldots \oplus U_r, \mbox{where $0 \neq U_k \subset W_{i_k}$ for some 
indices $1 \leq i_1 < \ldots < i_r \leq n$}
\]
with the subspaces $W_i$ from (\ref{eq:decompose-V2}). The indices $(i_1, \ldots, 
i_r)$ are called the {\em weights of $U$}.

\begin{prop} \label{prop:faithful}
Suppose that $H \backslash G/K$ is a minimal biquotient with non-central flats. Let $0 \neq U \subset 
\s^\perp$ be a $(\tau, \phi, \l)$-invariant subspace. Then $\l_s$ acts faithfully on $U$.
\end{prop}

\begin{proof}
Decompose $\l_s = \l_1 \oplus \l_2$ where $\l_2 := Ann(U) \cap \l_s \lhd \l_s$. We also 
decompose
\be \label{eq:decompose-H-1}
\H = \F \oplus (\H \cap \l_1) \oplus (\H \cap \l_2) \oplus \{v + \al(v) \mid v \in \Delta\}
\ee
for some linear isomorphism $\al: \Delta \ra \Delta'$ with $\Delta \subset \l_1 \cap 
(\H \cap \l_1)^\perp$ and $\Delta' \subset \l_2 \cap (\H \cap \l_2)^\perp$. Moreover, since 
$U$ and hence $\l_i$ are $\tau$-invariant, we have $\Delta, \Delta' \subset \k^\perp$.

Let $v \in \Delta$. Then $v_1 := v + \al(v) \in \H$, whereas $v_2 := v - (\al^t)^{-1} v \in \V \cap 
\k^\perp$, and $ad_{v_i}|_{\s^\perp} = J A_i$ for symmetric maps $A_i: V_2 \ra V_2$ for 
$i = 1,2$.  Since $\al(v), (\al^t)^{-1} (v) \in \Delta' \subset Ann(U)$, it follows that $(A_2)|_U = 
(A_1)|_U$.

But now, Lemma~\ref{lem:eigenvalues}.1 implies that $[P_0^{-1} A_1, \phi] = 0$, hence 
$[P_0^{-1} A_2, \phi]|_U = 0$ so that we can find a basis of $U$ consisting of common 
eigenvectors of $\phi$ and $P_0^{-1} A_2$. Therefore, Lemma~\ref{lem:eigenvalues}.2 
implies that $A_2(U) = 0$, so that $ad_v(U) = ad_{v_2}(U) = 0$ which implies that $v \in 
\Delta \cap \l_2 = 0$. Thus, $\Delta = \Delta' = 0$, i.e. (\ref{eq:decompose-H-1}) becomes
\[
\H = \F \oplus (\H \cap \l_1) \oplus (\H \cap \l_2).
\]
On the other hand,
\[
Q([\l_2, \s^\perp], U) = Q(\s^\perp, \underbrace{[\l_2, U]}_{=0}) = 0,
\]
hence $[\l_2, \s^\perp] \subset U^\perp \subsetneq \s^\perp$, so that 
Proposition~\ref{prop:l_s preserves s^perp} implies that $\l_2 = 0$ which shows the claim.
\end{proof}

Now define the orthogonal projections $\pi_k: \s^\perp \ra W_k$ w.r.t. the splitting 
(\ref{eq:decompose-V2}), and the maps 
\[
\phi_{ij} := \pi_i \circ \phi \circ \pi_j.
\]
Since $\phi$ is symmetric, it follows that $(\phi_{ij})^t = \phi_{ji}$. Also, 
$\phi = \sum_{i,j} \phi_{ij}$. For a tuple $\sigma = (i_1, \ldots, i_n)$ we define 
\[
\phi_\sigma =: \phi_{i_1 i_2} \circ \phi_{i_2 i_3} \circ \ldots \circ \phi_{i_{n-1}i_n},
\]
and we say that $\sigma$ is a tuple {\em from $i_1$ to $i_n$}. If $i_1 = i_n$ then we 
call $\sigma$ an {\em ($i_1$-based) loop}. Note that $(\phi_\sigma)^t = 
\phi_{\sigma^t}$, where $\sigma^t := (i_n, \ldots, i_1)$.

\begin{lem} \label{lem:m-phi-commute}
Let $\sigma = (i_1, \ldots, i_n)$ be a tuple from $i_1$ to $i_n$ and let  $v \in \H$. Then
\[
ad_v \circ \phi_\sigma = \frac {\la_{i_1}}{\la_{i_n}}  \phi_\sigma \circ ad_v 
\]
In particular, $\H$ and hence $\l$ commutes with $\phi_\sigma$ if $\sigma$ is a loop.
\end{lem}

\begin{proof} 
Evidently, it suffices to show the lemma for tuples $\sigma = (i, j)$. Let $A = A_v : V_2 
\ra V_2$ be the symmetric linear map such that $ad_v = JA$. By 
Lemma~\ref{lem:eigenvalues}, $[P_0^{-1} A, \phi] = 0$. Now let $A_i := A \circ \pi_i = 
\pi_i \circ A$. Then $P_0^{-1} A = \sum_i 1/\la_i A_i$, hence we have
\[
\ba{lll}
0 & = & \pi_i [P_0^{-1} A, \phi] \pi_j\\
& = & P_0^{-1} A_i \phi \pi_j - \pi_i \phi P_0^{-1} A_j\\
& = & 1/\la_i \ A_i \phi_{ij} - 1/\la_j\ \phi_{ij} A_j\\
& = & 1/\la_i\ A \phi_{ij} - 1/\la_j\ \phi_{ij} A,
\ea
\]
and from here the claim follows.
\end{proof}

A $(\tau, \phi, \l)$-invariant subspace $U \subset V_2$ is called {\em $(\tau, \phi, 
\l)$-irreducible} if it has no non-trivial $(\tau, \phi, \l)$-invariant subspace. If $U \subset 
\s^\perp$ is $(\tau, \phi, \l)$-invariant, then so is $\s^\perp \cap U^\perp$, hence 
$\s^\perp$ is the orthogonal sum of $(\tau, \phi, \l)$-irreducible subspaces.

\begin{df} \label{df:partial-order}
Let ${\cal I}$ be the set of all non-zero $\tau$-invariant ideals $\l' \lhd \l_s$ for which 
$(\l', \tau)$ is an irreducible symmetric pair. On ${\cal I}$, define the partial ordering 
$\less$ by
\[
\l_1 \less \l_2 \Longleftrightarrow \left\{\ba{l} \mbox{for all $v \in \H$, we have $pr_{\l_1}(v) 
\neq 0$ iff  $pr_{\l_2}(v) \neq 0$,}\\ \\ \mbox{and in this case, } ||pr_{\l_1}(v)||_{B_{\l_1}} < 
||pr_{\l_2}(v)||_{B_{\l_2}}, \ea \right.
\]
where for a compact semi-simple Lie algebra $\l'$, $||\cdot||_{B_{\l'}}$ denotes the norm 
w.r.t. the negative of the Killing form.
\end{df}

\noindent
Note that we have the decomposition $\l_s := \bigoplus_{\l' \in {\cal I}} \l'$.
 
\begin{prop} \label{prop:Max-Min}
Suppose that $H \backslash G/K$ is a minimal biquotient with non-central flats. Let $0 \neq U = U_1 
\oplus \ldots \oplus U_r \subset \s^\perp$ be a $(\tau, \phi, \l)$-irreducible subspace with 
weigths $i_1 < \ldots < i_r$, and let $\l_{\max}, \l_{\min} \lhd \l_s$ be maximal and minimal 
ideals w.r.t. the partial ordering $\less$. Then
\bi
\item
$[\l_{\min}, U_1] = U_1$, and $[\l_{\min}, U_k] = 0$ for all $k > 1$,
\item
$[\l_{\max}, U_r] = U_r$, and $[\l_{\max}, U_k] = 0$ for all $k < r$.
\ei
\end{prop}

\begin{proof}
Suppose there is a subspace $0 \neq U_1' \subset U_1$ which is invariant under $\tau$, 
$\l$ and under $\phi_\sigma$ for all loops $\sigma$. Let $U_k' := span \{\phi_\sigma(U_1') 
\mid \mbox{$\sigma$ a tuple from $i_1$ to $i_k$}\} \subset U_k$. Then $U_k'$ is 
$\tau$-invariant, and $\phi_{ij}(U_j') \subset U_i'$ by the definition of $\phi_\sigma$. Also, 
since for all $v \in \H$, $\phi_\sigma ad_v$ and $ad_v \phi_\sigma$ are multiples of each 
other by Lemma~\ref{lem:m-phi-commute}, $U_k'$ is also $\l$-invariant, so that 
$\bigoplus_k U_k' \subset U$ is $(\tau, \phi, \l)$-invariant. But since we assume that $U$ 
is $(\tau, \phi, \l)$-irreducible, this means that this direct sum is all of $U$ and hence $U_k = 
U_k'$ for all $k$.

Decompose $U_1 := \sum_\al U_1^\al \ot \R^{n_\al}$ into $\tau$- and $\l$-invariant 
subspaces, such that $U_1^\al$ and $U_1^\beta$ are inequivalent $\tau$-equivariant 
irreducible representations of $\l$ for $\al \neq \beta$. Since $\tau$ and $\l$ commute 
with $\phi_\sigma$ for all loops $\sigma$, it follows that the summands  $U_1^\al \ot 
\R^{n_\al}$ are invariant under $\tau$, $\l$ and these $\phi_\sigma$, so that the preceding 
paragraph implies that $U_1$ contains only one such summand. Arguing likewise for all 
$U_k$, we conclude $U_k = \tilde U_k \ot \R^{n_k}$ for some irreducible $\tau$-equivariant representation $\tilde U_k$ of $\l$.

For each $k$, pick a tuple $\sigma_k$ from $i_1$ to $i_k$ such that $\phi_{\sigma_k}(U_1) 
\neq 0$. Then there is an $e_1 \in \R^{n_1}$ such that $(\phi_{\sigma_k})|_{\tilde U_1 \ot e_1} 
\neq 0$ for all $k$. Note that $\ker (\phi_{\sigma_k}) = \ker \left((\phi_{\sigma_k})^t 
\phi_{\sigma_k}\right) = \ker \phi_{\sigma_k * \sigma_k^t}$, and since $\sigma_k * 
\sigma_k^t$ is a loop, it follows that $\ker (\phi_{\sigma_k}) \subsetneq \tilde U_1 \ot e_1$ is 
a $\tau$- and $\l$-invariant proper subspace which must vanish. That is, 
$\phi_{\sigma_k}|_{\tilde U_1 \ot e_1}$ is injective.

By Lemma~\ref{lem:m-phi-commute}, $ad_v \circ \phi_{\sigma_k}$ and $\phi_{\sigma_k} \circ 
ad_v$ are multiples of each other for all $v \in \H$, hence $\phi_{\sigma_k}(\tilde U_1 \ot e_1) 
\subset U_k$ is $\tau$- and $\l$-invariant and irreducible, that is, $\phi_{\sigma_k}(\tilde U_1 
\ot e_1) = \tilde U_k \ot e_k \subset U_k$ for some $0 \neq e_k \in \R^{n_k}$. Note that $e_k$ 
is well-defined only up to scalar multiples. However, once we fixed these elements $e_k$, 
we may identify $\tilde U_k \cong \tilde U_k \ot e_k$ for $k = 1, \ldots, r$ and hence the adjoint representation of $\l_s$ on $U$ gives rise to a homomorphism
\be \label{eq:representation}
ad_\l: \l_s \longrightarrow \bigoplus_{k=1}^r End(\tilde U_k),
\ee
and we may regard $\phi_{\sigma_k}$ as an isomorphim $\phi_{\sigma_k}: 
\tilde U_1 \ra \tilde U_k$.

We decompose the projection of $ad_\l(\l_s)$ to $End(\tilde U_1)$ into irreducible 
$\tau$-invariant ideals $\uu_1^1 \oplus \ldots \oplus \uu_1^m$, and define $\uu_k^\nu := 
Ad_{\phi_{\sigma_k}}(\uu_1^\nu) \subset End(\tilde U_k)$ for $k = 1, \ldots, r$ and $\nu = 1, 
\ldots, m$. Once again by Lemma~\ref{lem:m-phi-commute},
\be \label{eq:commuting}
Ad_{\phi_{\sigma_k}}^{-1} \left((ad_v)|_{\tilde U_k}\right) = (\phi_{\sigma_k}^{-1} \circ ad_v 
\circ \phi_{\sigma_k})|_{\tilde U_1} = \frac{\la_{i_1}}{\la_{i_k}} (ad_v)|_{\tilde U_1},
\ee
so that the projection of $ad_\l(\l_s)$ to $End(\tilde U_k)$ equals $\uu_k^1 \oplus \ldots 
\oplus \uu_k^m$. Therefore, we can refine the faithful representation (\ref{eq:representation}) to
\[
ad_\l: \l_s \longrightarrow \bigoplus_{k=1}^r \bigoplus_{\nu=1}^m \uu_k^\nu,
\]
where the projection to each summand is surjective. Since $\l_s$ acts faithfully on $U$ by 
Proposition~\ref{prop:faithful}, it follows that $ad_\l^{-1}(\uu_k^\nu) \lhd \l_s$ is an ideal which is 
$\tau$-equivariantly isomorphic to $\uu_k^\nu$ and hence irreducible, i.e., $ad_\l^{-1}(\uu_k^\nu) 
\in {\cal I}$ with the set ${\cal I}$ from Definition~\ref{df:partial-order}. Thus, we get a map
\[
I_0: \{1, \ldots, r\} \times \{1, \ldots, m\} \longrightarrow {\cal I},\ \ I_0(k, \nu) := 
ad_\l^{-1}(\uu_k^\nu) \in {\cal I}.
\]

Let $v \in \H$, and denote by $u_k^\nu \in \uu_k^\nu \subset End(\tilde U_k)$ the projection 
of $(ad_v)|_{\tilde U_k}$. Then (\ref{eq:commuting}) implies that
\be \label{eq:commuting2}
Ad_{\phi_{\sigma_k}}^{-1}(u_k^\nu) = \frac{\la_{i_1}}{\la_{i_k}} u_1^\nu.
\ee
The Killing form is preserved under Lie algebra isomorphisms, and since both $ad_\l: 
I_0(k, \nu) \ra \uu_k^\nu$ and $Ad_{\phi_{\sigma_k}}^{-1}: \uu_k^\nu \ra \uu_1^\nu$ are 
isomorphisms, we have
\[ 
||pr_{I_0(k, \nu)}(v)||_{B_{I_0(k, \nu)}} = ||u_k^\nu||_{B_{\uu_k^\nu}} = 
||Ad_{\phi_{\sigma_k}}^{-1} (u_k^\nu)||_{B_{\uu_1^\nu}} = \frac{\la_{i_1}}{\la_{i_k}} 
||u_1^\nu||_{B_{\uu_1^\nu}},
\]
where we used (\ref{eq:commuting2}) in the last equation. From this, we deduce 
immediately that
\[
||pr_{I_0(k, \nu)}(v)||_{B_{I_0(k, \nu)}} = \frac{\la_{i_{k'}}}{\la_{i_k}} ||pr_{I_0(k', 
\nu)}(v)||_{B_{I_0(k', \nu)}},
\]
and since for $k < k'$ we have $\la_{i_k} < \la_{i_{k'}}$, it follows that
\[
I_0(k, \nu) \less I_0(k', \nu) \ \mbox{for all $k < k'$ and all $\nu$}.
\]

Thus, if $\l_{\min} \in {\cal I}$ is minimal w.r.t. $\less$, then 
$pr_{\uu_k^\nu}(ad_\l(\l_{\min})) = 0$ for all $k > 1$, that is, $[\l_{\min}, \tilde U_k] = 0$ 
and therefore $[\l_{\min}, U_k] = 0$ for all $k > 1$. Since $\l_{\min} \lhd \l_s$ is a 
$\tau$-invariant ideal, it follows that $[\l_{\min}, \tilde U_1] \subset \tilde U_1$ is $\tau$- 
and $\l$-invariant, hence by irreducibility is either zero or all of $\tilde U_1$. In the first 
case, however, we would have $ad_{\l_{\min}}|_{s^\perp} = 0$, contradicting the 
faithfulness of the action of $\l_s$, hence $[\l_{\min}, \tilde U_1] = \tilde U_1$ and 
therefore $[\l_{\min}, U_1] = U_1$ as claimed.

The corresponding statements about $\l_{\max}$ follow analogously.
\end{proof}

\noindent Proposition~\ref{prop:Max-Min} is the last step which is needed to show the

\begin{thm} \label{thm:MainTheorem}
There are no normal biquotients with non-central flats.
\end{thm}

\begin{proof}
Suppose biquotients with non-central flats exist. Then there exists also a {\em minimal} 
biquotient with non-central flats $H \backslash G/K$. Let $\l_{\max}, \l_{\min} \lhd \l_s \subset \s$ and $0 \neq 
U = U_1 \oplus \ldots U_r \subset \s^\perp$ be as in Proposition~\ref{prop:Max-Min}. If 
$r = 1$, then any eigenvector of $\phi|_U$ would also be an eigenvector of $P_0$ as 
$U \subset W_{i_1}$. But this is impossible by Lemma~\ref{properties phi}.

Thus, $r > 1$ and hence $[\l_{\max}, [\l_{\min}, U]] = 0$ by Proposition~\ref{prop:Max-Min}. 
Since $\s^\perp$ is the direct sum of $(\tau, \phi, \l)$-irreducible subspaces, it follows that
\[
[\l_{\max}, [\l_{\min}, \s^\perp]] = 0.
\]

Also, $W_1 = \bigoplus \{U_1 \mid 0 \neq U_1 \oplus \ldots \oplus U_r \mbox{ is 
$(\tau, \phi, \l)$-irreducible with $1$ as a weight}\}$. Again by 
Proposition~\ref{prop:Max-Min}, this implies that
\[
[\l_{\min}, W_1] = W_1.
\]

By Proposition~\ref{prop:irreducible}, $G/K$ is a compact irreducible symmetric space, 
and if we apply Lemma~\ref{lem:composition=0} with $\l_1 := \l_{\min}$, $\l_2 := \l_{\max}$ 
and $W := W_1$, then we conclude $\l_{\max} = 0$ which is a contradiction.
\end{proof}


{\sc
\noindent
Fachbereich Mathematik, Universit\"at Dortmund, 44221 Dortmund, Germany
\

\noindent
Email:\, 
{\tt lschwach@math.uni-dortmund.de}
}

\end{document}